\documentclass{amsart}
\usepackage{amssymb}
\usepackage{graphicx}

\newtheorem{thm}{Theorem}
\newtheorem{cor}{Corollary}
\newtheorem{lem}{Lemma}

\newtheorem{rem}{Remark}

\newcommand{\D}{{\mathbb D}}

\newcommand{\real}{{\operatorname{Re}\,}}

\def\be{\begin{equation}}
\def\ee{\end{equation}}

\newcommand{\bee}{\begin{enumerate}}
\newcommand{\eee}{\end{enumerate}}

\newcommand{\blem}{\begin{lem}}
\newcommand{\elem}{\end{lem}}
\newcommand{\bthm}{\begin{thm}}
\newcommand{\ethm}{\end{thm}}
\newcommand{\bcor}{\begin{cor}}
\newcommand{\ecor}{\end{cor}}
\newcommand{\beg}{\begin{example}}
\newcommand{\eeg}{\end{example}}
\newcommand{\begs}{\begin{examples}}
\newcommand{\eegs}{\end{examples}}
\newcommand{\bdefe}{\begin{defin}}
\newcommand{\edefe}{\end{defin}}
\newcommand{\bprob}{\begin{prob}}
\newcommand{\eprob}{\end{prob}}
\newcommand{\bei}{\begin{itemize}}
\newcommand{\eei}{\end{itemize}}

\newcommand{\bcon}{\begin{conj}}
\newcommand{\econ}{\end{conj}}
\newcommand{\bcons}{\begin{conjs}}
\newcommand{\econs}{\end{conjs}}
\newcommand{\bprop}{\begin{propo}}
\newcommand{\eprop}{\end{propo}}
\newcommand{\br}{\begin{rem}}
\newcommand{\er}{\end{rem}}
\newcommand{\brs}{\begin{rems}}
\newcommand{\ers}{\end{rems}}
\newcommand{\bo}{\begin{obser}}
\newcommand{\eo}{\end{obser}}
\newcommand{\bos}{\begin{obsers}}
\newcommand{\eos}{\end{obsers}}
\newcommand{\bpf}{\begin{pf}}
\newcommand{\epf}{\end{pf}}
\newcommand{\ba}{\begin{array}}
\newcommand{\ea}{\end{array}}
\newcommand{\beq}{\begin{eqnarray}}
\newcommand{\beqq}{\begin{eqnarray*}}
\newcommand{\eeq}{\end{eqnarray}}
\newcommand{\eeqq}{\end{eqnarray*}}

\begin{document}
\bibliographystyle{amsplain}

\title[The third logarithmic coefficient for the class $\mathcal{S}$]{The third logarithmic coefficient for the class $\boldsymbol{\mathcal{S}}$}

\author[M. Obradovi\'{c}]{Milutin Obradovi\'{c}}
\address{Department of Mathematics,
Faculty of Civil Engineering, University of Belgrade,
Bulevar Kralja Aleksandra 73, 11000, Belgrade, Serbia}
\email{obrad@grf.bg.ac.rs}

\author[N. Tuneski]{Nikola Tuneski}
\address{Department of Mathematics and Informatics, Faculty of Mechanical Engineering, Ss. Cyril and Methodius
University in Skopje, Karpo\v{s} II b.b., 1000 Skopje, Republic of North Macedonia.}
\email{nikola.tuneski@mf.edu.mk}

%\author{M. Obradovi\'{c}}
%\address{M. Obradovi\'{c},
%Department of Mathematics,
%Faculty of Civil Engineering, University of Belgrade,
%Bulevar Kralja Aleksandra 73, 11000
%\email{obrad@grf.bg.ac.rs}

%\author{S. Ponnusamy${}^{~\mathbf{*}}$}
%\address{S. Ponnusamy, Department of Mathematics,
%Indian Institute of Technology Madras, Chennai--600 036, India.}
%\email{samy@iitm.ac.in}
%\author{N. Tuneski}
%\address{N. Tuneski, St. Cyril and Methodius University, Faculty of Mechanical
%Engineering, Karpo\v s II b.b., 1000 Skopje, R. Macedonia.}
%\email{nikolat@mf.edu.mk}

%\date{\today  %June. 30, 09
%;  File: }

%\begin{abstract}

%\end{abstract}

%\thanks{The work of the first author was supported by MNZZS Grant, No. ON174017, Serbia. The research of the second
%author was supported by National Board for Higher
%Mathematics, India.}

\begin{abstract}
In this paper we give an upper bound  of the third logarithmic coefficient for the class  $\mathcal{S}$ of univalent functions in the unit disc.
\end{abstract}

\maketitle
%\pagestyle{myheadings}
%\markboth{M. Obradovi\'{c} and  S. Ponnusamy }{}
%\cc
%\section{Introduction}

\medskip

%\section{Introduction and preliminaries}

Let $\mathcal{A}$ be the class of functions $f$ that are analytic  in the open unit disc $\D=\{z:|z|<1\}$ of the form
\be\label{eq 1}
f(z)=z+a_2z^2+a_3z^3+\cdots,
\ee
and let $\mathcal{S}$ be its subclass consisting of functions that are univalent in the unit disc $\D$.

\medskip

The logarithmic coefficients of the function $f$ given by \eqref{eq 1} are defined in $\D $ by
\be\label{eq 2}
\log\frac{f(z)}{z}=2\sum_{n=1}^{\infty}\gamma_{n}z^{n}.
\ee
By using  \eqref{eq 1}, after differentiation and comparing the coefficients, we can obtain that $\gamma_1=\frac12a_2$, $\gamma_2=\frac12\left(a_3-\frac12a_2^2\right)$ and
\be\label{eq 3}
\gamma_{3}=\frac{1}{2}\left(a_{4}-a_{2}a_{3}+\frac{1}{3}a_{2}^{3}\right).
\ee

Very little is known about the estimates of the modulus of the logarithmic coefficients for the whole class $\mathcal{S}$ of normalized of univalent functions. The Koebe function $k(z)=\frac{z}{(1-z)^2}=\sum_{n=1}^\infty z^n$ with $\gamma_n=\frac1n$ being extremal in majority estimates over the class $\mathcal{S}$ inspires a conjecture that $|\gamma_n|\le\frac1n$ for $n=1,2,\ldots$ and $f\in\mathcal{S}$. Apparently, this is true only for the class of starlike functions (\cite{Thomas-2}), but not for the class $\mathcal{S}$ in general (\cite[Thgeorem 8.4, p.242]{duren}). Sharp estimates for the class $\mathcal{S}$ are known only for the first two coefficients, $|\gamma_1|\le1$ and $|\gamma_2|\le \frac12+\frac1e$.

\medskip

In this paper we give an upper bound of $|\gamma_{3}|$ for the class $\mathcal{S}$.

\medskip

It is worth mentioning that the problem of estimating the modulus of the first three logarithmic coefficients is widely studied for the subclasses of $\mathcal{S}$ and in some cases sharp bounds are obtained. Namely,
sharp estimates for the class of strongly starlike functions of certain order and $\gamma$-starlike functions are given in \cite{Thomas-2} and \cite{darus}, respectively, while non-sharp estimates for the class of Bazilevic, close-to-convex and different subclasses of close-to-convex functions are given in \cite{deng}, \cite{ali} and \cite{Thomas}, respectively.

%
%\medskip
%
%In 2016 Thomas (\cite{Thomas}) found sharp estimate  $|\gamma_{3}|\leq\frac{7}{12}=0.5833\ldots$ for certain subclass of close-to-convex functions.
%Also, recently (\cite{Cho}) the authors considered different subclasses of close-to-convex
%functions and found appropriate sharp results
%
%\medskip

%\section{Main result}

As announced before, here is an estimate of the modulus of the third logarithmic coefficient for the whole class of univalent functions.

\medskip

\bthm\label{20-th 1} For the class $\mathcal{S}$ we have
$$
|\gamma_{3}|\leq \frac{\sqrt{133}}{15} = 0.7688\ldots.
$$
\ethm

\medskip

\begin{proof}
In the proof of this theorem we will use mainly the notations and results given in the book of N. A. Lebedev (\cite{Lebedev}).

\medskip

Let $f \in \mathcal{S}$ and let
\[
\log\frac{f(t)-f(z)}{t-z}=\sum_{p,q=0}^{\infty}\omega_{p,q}t^{p}z^{q},
\]
where $\omega_{p,q}$ are called Grunsky's coefficients with property $\omega_{p,q}=\omega_{q,p}$.
For those coefficients we have the next Grunsky's inequality (\cite{duren,Lebedev}):
\be\label{eq 4}
\sum_{q=1}^{\infty}q \left|\sum_{p=1}^{\infty}\omega_{p,q}x_{p}\right|^{2}\leq \sum_{p=1}^{\infty}\frac{|x_{p}|^{2}}{p},
\ee
where $x_{p}$ are arbitrary complex numbers such that last series converges.

\medskip

Further, it is well-known that if $f$ given by \eqref{eq 1}
belongs to $\mathcal{S}$, then also
\be\label{eq 5}
f_{2}(z)=\sqrt{f(z^{2})}=z +c_{3}z^3+c_{5}z^{5}+\cdots
\ee
belongs to the class $\mathcal{S}$. Then for the function $f_{2}$ we have the appropriate Grunsky's
coefficients of the form $\omega_{2p-1,2q-1}^{(2)}$ and the inequality \eqref{eq 4} has the form
\be\label{eq 6}
\sum_{q=1}^{\infty}(2q-1) \left|\sum_{p=1}^{\infty}\omega_{2p-1,2q-1}^{(2)}x_{2p-1}\right|^{2}\leq \sum_{p=1}^{\infty}\frac{|x_{2p-1}|^{2}}{2p-1}.
\ee
As it has been shown in \cite[p.57]{Lebedev}, if $f$ is given by \eqref{eq 1} then the coefficients $a_{2}, a_{3}, a_{4}$
are expressed by Grunsky's coefficients  $\omega_{2p-1,2q-1}^{(2)}$ of the function $f_{2}$ given by
\eqref{eq 5} in the following way (in the next text we omit upper index 2 in $\omega_{2p-1,2q-1}^{(2)}$):

\be\label{eq 7}
\begin{split}
a_{2}&=2\omega _{11},\\
a_{3}&=2\omega_{13}+3\omega_{11}^{2}, \\
a_{4}&=2\omega_{33}+8\omega_{11}\omega_{13}+\frac{10}{3}\omega_{11}^{3}.
\end{split}
\ee
Now, from \eqref{eq 3} and \eqref{eq 7} we have
$$
\gamma_{3}=\omega_{33}+2\omega_{11}\omega_{13}
$$

\medskip

On the other hand, from \eqref{eq 7} for $x_{2p-1}=0$, $p=3,4,\ldots$ we have
\be\label{eq 9}
|\omega_{11}x_{1}+\omega_{31}x_{3}|^{2}+3|\omega_{13}x_{1}+\omega_{33}x_{3}|^{2}
\leq |x_{1}|^{2}+\frac{|x_{3}|^{2}}{3}.
\ee
From \eqref{eq 9} for $x_{1}=2\omega_{11}$, $x_{3}=1$ and since $\omega_{31}=\omega_{13}$, we have
$$|2\omega_{11}^2+\omega_{13}|^{2} +3| \gamma_3|^2 \leq 4|\omega_{11}|^2+\frac13,$$
and from here
\[
\begin{split}
 |\gamma_3|^2
 & \le \frac19 +\frac43|\omega_{11}|^2-\frac13|2\omega_{11}^2+\omega_{13}|^2 \\
 & = \frac19 + \frac43|\omega_{11}|^2-\frac13\left( 4|\omega_{11}|^4 + |\omega_{13}|^2 +4\real\left\{ \omega_{13}\overline{\omega_{11}}^2 \right\} \right) \\
 & = \frac19 + \frac43|\omega_{11}|^2-\frac43 |\omega_{11}|^4 -\frac13 |\omega_{13}|^2 -\frac43 \real\left\{ \omega_{13}\overline{\omega_{11}}^2 \right\}.
\end{split}
\]
Using the fact that
\[  -|\omega_{13}|^2 \le -\left|\real\{\omega_{13}\}\right|^2 = -\left(\real\{\omega_{13}\}\right)^2, \]
we obtain
\[
 |\gamma_3|^2 \le \frac19 + \frac43|\omega_{11}|^2-\frac43 |\omega_{11}|^4 - \frac13\left(\real\{\omega_{13}\}\right)^2 -\frac43 \real\left\{ \omega_{13}\overline{\omega_{11}}^2 \right\}.
\]

\medskip

Next, without loss of generality using suitable rotation of $f$ we can assume that $0\le a_2 \le2$  and $a_2=2\omega_{11}$ receive that $0\le \omega_{11}\le1$.
So, let put $\omega_{11}=a$, $0\le a\le1$, and continue analysing
\begin{equation}\label{eq10}
 |\gamma_3|^2 \le \frac19 + \frac43a^2-\frac43 a^4 -\frac13\left(\real\{\omega_{13}\}\right)^2 -\frac43a^2 \real\left\{ \omega_{13}\right\}.
\end{equation}

\medskip

It is a classical result that for the class $\mathcal{S}$ we have $|a_3-a_2^2|\le1$ (see \cite[p.5]{DTV-book}), which is by \eqref{eq 7} equivalent with
\[ |2\omega_{13}-\omega_{11}^2| \le 1. \]
From here,
\[ -1 \le \real \{ 2\omega_{13}-\omega_{11}^2 \} \le 1, \]
i.e.,
\begin{equation}\label{eq11}
-\frac12(1-a^2) \le \real \{ \omega_{13}\} \le \frac12(1+a^2).
\end{equation}

\medskip

If we put $x_1=1$ and $x_3=0$ in \eqref{eq 9}, then we get
\[ |\omega_{11}|^2+3|\omega_{13}|^2\le 1, \]
which implies
\[ |\omega_{13}|\le \frac{1}{\sqrt3}\sqrt{1-|\omega_{11}|^2 }= \frac{1}{\sqrt3} \sqrt{1-a^2}.  \]
Combining this with \eqref{eq11}, we receive
\[ -\frac12(1-a^2) \le \real \{\omega_{13}\} \le \frac{1}{\sqrt{3}}\sqrt{1-a^2} \]
(because $-\frac{1}{2}(1-a^2) \ge -\frac{1}{\sqrt3}\sqrt{1-a^2} $).

\medskip

By using \eqref{eq10}, \eqref{eq11} and the notation $t=\real \{ \omega_{13}\}$ we obtain
\[
|\gamma_3|^2 \le \frac19 + \frac43a^2-\frac43 a^4 -\frac13t^2 -\frac43a^2 t :\equiv \psi(a,t) = \frac19+\frac13\varphi(a,t),
\]
where $0\le a\le1$,  $-\frac12(1-a^2) \le t \le \frac{1}{\sqrt{3}}\sqrt{1-a^2}$ and $\varphi(a,t)=4a^2-4a^4-t^2-4a^2t$.

\medskip

It remains to show that the maximal value of the function $\psi(a,t)$ over the region $\Omega=[0,1]\times[-\frac12(1-a^2),\frac{1}{\sqrt{3}}\sqrt{1-a^2}]$ equals $\left(\frac{\sqrt{133}}{15}\right)^2 = \frac{133}{225}$, or equivalently that $\varphi(a,t)$ has maximal value $\frac{36}{25}$ on the same region.

\medskip

Indeed, the system of equations
\[
\left\{
\begin{array}{l}
\varphi'_a(a,t)=8a-16a^3-8at=0\\
\varphi'_t(a,t)=-4a^2-2t=0
\end{array}
\right.
\]
has unique real solution $a=t=0$ with $\varphi(0,0)=0$, while on the edges of the region $\Omega$ we have the following:
\begin{itemize}
  \item[-] for $a=0$ we have that the function  $\varphi(0,t)=-t^2$ on the interval $-\frac12\le t\le\frac{1}{\sqrt3}$ attains maximal value $\varphi(0,0)=0$;
  \item[-] when $a=1$, $t$ can take single value, $t=0$, and in that case $\varphi(1,0)=0$;
  \item[-] for $t=-\frac12(1-a^2)$, the function $\varphi\left(a,-\frac12(1-a^2)\right)= -\frac14 (a^2-1)\left(a^2-\frac{1}{25}\right)$ is with maximal value $\frac{36}{25}$ on the interval $0\le a\le1$ attained for $a=\frac{\sqrt{13}}{5}$;
  \item[-] for $t=\frac{1}{\sqrt{3}}\sqrt{1-a^2}$, the values of the function
  \[
  \begin{split}
  \varphi\left(a,\frac{1}{\sqrt{3}}\sqrt{1-a^2}\right)
  &= \frac13(-12a^4+13a^2-1)-\frac{4a^2}{\sqrt3}\sqrt{1-a^2} \\
  &\le \frac13(-12a^4+13a^2-1)<\frac{36}{25}.
  \end{split}
  \]
  on the interval $0\le a\le1$ are smaller than $\frac{36}{25}$.
\end{itemize}

\medskip

This completes the proof.
\end{proof}

\end{document}